\documentclass[11pt]{amsart}
\usepackage{fullpage}
\usepackage{amsfonts}
\usepackage{latexsym}
\usepackage{amscd}
\usepackage{amsmath}
\usepackage{amssymb}
\usepackage{graphicx}
\usepackage[mathscr]{eucal}
\usepackage{xypic}
\usepackage{color}

\newtheorem{teo}{Theorem}[section]
\newtheorem{coro}[teo]{Corollary}
\newtheorem{lema}[teo]{Lemma}

\newtheorem{ejem}[teo]{Example}

\newtheorem{question}[teo]{Question}

\newtheorem{problema}[teo]{Problem}

\def\nat{\mathbb{N}}

\def\nin{\not\in}

\def\w1{\omega_1}

\begin{document}

\title{On the continuity of the elements of the Ellis semigroup and other properties}

\author{S. Garc\'{\i}a-Ferreira}
\address{Centro de Ciencias Matem\'aticas, Universidad Nacional
Aut\'onoma de M\'exico, Apartado Postal 61-3, Santa Mar\'{\i}a,
58089, Morelia, Michoac\'an, M\'exico}
\email{sgarcia@matmor.unam.mx}

\author{Y. Rodr\'iguez-L\'opez}
\address{Departamento de Matem\'atica, Instituto de Ci\^encias Exatas, Universidade Federal de Minas Gerais, Av. Ant\^onio Carlos, 6627, Belo Horizonte, Brazil}
\email{yackelin@mat.ufmg.br}

\author{C. Uzc\'ategui}
\address{Escuela de Matem\'aticas, Facultad de Ciencias, Universidad Industrial de
Santander, Ciudad Universitaria, Carrera 27 Calle 9, Bucaramanga,
Santander, A.A. 678, Colombia and
Departamento de Matem\'{a}ticas,
Facultad de Ciencias, Universidad de los Andes, M\'erida 5101,
Venezuela}
\email{cuzcatea@saber.uis.edu.co}

\thanks{Research of the first-named author was supported by the PAPIIT grant no. IN-105318.  The third author was supported by the  VIE-UIS  grant \#2422.}



\dedicatory{}


\begin{abstract} We consider discrete dynamical systems whose phase spaces are compact metrizable countable spaces.  In the first part of the article,  we study  some properties that guarantee the continuity of all  functions of the corresponding Ellis semigroup. For instance, if every accumulation point of $X$  is fixed,  we give a necessary and sufficient condition on a point  $a\in X'$ in order that all functions of the Ellis semigroup $E(X,f)$ be continuous at the given point $a$. In the second part,  we consider  transitive dynamical systems. We show that if $(X,f)$ is a transitive dynamical system and either every function of $E(X,f)$ is continuous or  $|\omega_f(x)|=1$ for each accumulation point $x$ of $X$, then  $E(X,f)$ is homeomorphic to $X$.  Several examples are given to illustrate our results. 
\end{abstract}

\maketitle

\section{Introduction and preliminaries}

 A dynamical systems $(X,f)$ will consist of a compact metric space  $X$
 and a continuous function $f:X\rightarrow X$ (usually, this kind of dynamical systems are called discrete dynamical systems and we short the name for convenience).  The {\it orbit} of point $x \in X$ is the set $\mathcal O_f(x) := \{ f^n(x) : n \in \mathbb{N} \}$.   The system $(X,f)$ is  {\it transitive} if there is a point with dense orbit. A point $x\in X$ is called
{\it periodic } if there is  $n\geq 1$ such that $f^n(x)=x$, and its {\it period} is  $\min\{n\in\mathbb N:f^n(x)=x\}$.
The symbol  $P_f$ stands for  the set of all periods of the periodic points of a dynamical systems $(X,f)$. A point $x$ is called {\it eventually periodic} if its orbit is finite. The $\omega-${\it limit set} of $x\in X$, denoted by $\omega_f(x)$, is the set of points $y\in X$ for which there exists an increasing sequence $(n_k)_{k\in\mathbb N}$ such that $f^{n_k}(x)\rightarrow y$. Observe that for each $y\in\mathcal O_f(x)$, we have that
$\omega_f(y)=\omega_f(x)$. If $\mathcal{O}_f(y)$ contains a periodic point $x$, then $\omega_f(y)=\mathcal{O}_f(x)$.
For a space $X$ we denote by $\mathcal{N}(x)$ the set of all  neighborhoods of $x\in X$, and the set of all accumulation points of $X$ will be simply denoted by $X'$. The Stone-\v{C}ech compactification $\beta(\mathbb N)$ of $\mathbb N$ with the discrete topology will be identified
with the set of ultrafilters over $\mathbb N$. Its remainder is denoted by ${\mathbb N}^*= \beta({\mathbb N})\setminus \mathbb{N}$ is  the
set of all free ultrafilters on $\mathbb N$, where, as usual, each natural number $n$ is identified with  the fixed ultrafilter
consisting of all subsets of $\mathbb N$ containing $n$.

\medskip

Since our phase spaces are compact metric countable spaces, we  remind the reader the classical  result from  \cite{ms} which asserts that  every  compact metric countable space  is homeomorphic to a countable sucesor ordinal with the order topology. In this context, some of the most attractive phase space has the form
$\omega^\alpha +1$, where $\alpha \geq 1$ is a countable ordinal.

\medskip

The {\it Ellis semigroup}  of a dynamical system $(X,f)$, denoted $E(X,f)$, is defined as  the pointwise closure of $\{f^n:\; n\in \mathbb N\}$ in the compact space $X^X$ with composition of functions as its algebraic operation.  The Ellis semigroup is equipped with the topology inhered from the product space $X^X.$ This semigroup was introduced by R. Ellis in \cite{ell}, and it is a very important tool in the  study of the topological behavior of a dynamical systems.  The article \cite{glasner} offers an excellent survey concerning applications of the Ellis semigroup. In the paper \cite{GarciaSanchis}, the authors initiated the study of  the continuity and discontinuity of the elements of $E(X,f) \setminus  \{ f^n : n \in \mathbb{N} \}$. For instance, they point out that  if $X$ is a convergent sequence with its limit point, then all the elements of $E(X,f)$ are either continuous or discontinuous (this result was later  improved in \cite{gru}). In a different context,  P. Szuca \cite{Szuca} showed that if $X = [0,1]$, the function $f: [0,1] \to [0,1]$ is onto and $f^p$ is continuous for some $p \in \mathbb{N}^*$, then all the elements of $E([0,1],f)$ are continuous. Using the Cantor set as a phase space and a generalized ship maps, the continuity and discontinuity  of the elements of the Ellis semigroup where studied  in \cite{gf}. The main  tool that have been used in all these investigations   is the combinatorial properties of the ultrafilters on $\mathbb{N}$. Certainly, the Ellis semigroup can be described  in terms of the notion of  convergence with respect to an ultrafilter. Indeed, given $p\in \mathbb{N}^*$ and a sequence $(x_n)_{n\in\mathbb{N}}$ in a space $X$, we say that a point $x\in X$ is the $p-${\it limit
point} of the sequence, in symbols $x =p-\lim_{n\rightarrow \infty}x_n$, if for every neighborhood $V$ of $x$, $\{n\in\mathbb{N}: f^n(x)\in V\} \in p$.
Observe that a point $x\in X$ is an accumulation point of a countable set $\{x_n:\,n\in\mathbb{N}\}$ of $X$ iff there is $p\in \mathbb{N}^*$ such that $x = p-\lim_{n\rightarrow \infty}x_n$.

\medskip

The notion of a $p-$limit point has been used in several branches of mathematics (see for
instance \cite{Be} and \cite[p. 179]{fu}). A. Blass \cite{Bla} and N. Hindman \cite{hi} formally established
the connection between  ``the iteration in to\-po\-lo\-gi\-cal dynamics'' and
``the convergence with respect to an ultrafilter'' by considering a more
general iteration of the function $f$ as follows: Let $X$ be
 a compact  space and $f : X \rightarrow  X$ a continuous function. For $p\in\mathbb{N}^*$, the
$p-${\it iterate} of $f$ is the function $f^p: X\rightarrow X$ defined by \[ f^p(x) = p-\lim_{n\rightarrow \infty} f^n(x),
\]
for each $x\in X$. The description of the Ellis semigroup and its operation in terms of the $p-$iterates are the following  (see \cite{Bla}, \cite{hi}):
\[
\begin{array}{rcl}
E(X,f) & = &\{f^p: p\in \beta\mathbb N \}  \\
\\
f^p\circ f^q & =& f^{q+p}\;\; \mbox{for each $p,q\in\beta\mathbb
N$}.
\end{array}
\]
We will use the following notation
\[
E(X,f)^* :=E(X,f)\setminus\{f^n:n\in\mathbb N\}.
\]
Besides, we have that $\omega_f(x)=\{f^p(x):p\in\mathbb N^*\}$ for each $x \in X$.

\medskip

Recently (for instance see \cite{gru} and \cite{gru2016}), we have investigated the structure of the Ellis semigroup of a dynamical system  and the topological properties of some of its elements.
Our main purpose moves in two directions: The first one concerns with the continuity and discontinuity of the $p$-iterates which is dealt in the third section, and the second one concerns about  a general question stated in \cite{gru2016}:

\begin{question}  Given two compact metric countable spaces $X$ and $Y,$ is there a continuous function $f : X \to X$ such that
$E(X,f)$  is homeomorphic to $Y$?
\end{question}

\noindent We provide a partial answer to this question in the forth section.

\section{Basic Properties}

In this section, we state several useful results that were proved in the articles \cite{gru} and  \cite{gru2016}. Our first lemma is
precisely  Lemma 2.1 from \cite{gru2016}.

\begin{lema}\label{aisladoperiodico} Let $(X,f)$ be a dynamical system and $x \in X$.
\begin{itemize}
\item[$(i)$] Assume that $x$ is periodic with period $n$ and let $l < n$. Then,   $p \in \big(n
\mathbb{N} + l\big)^*$ iff  $f^p(x) = f^l(x)$.

\item[$(ii)$]  Suppose that  $x$  is eventually periodic and that $m \in
\nat$  is the smallest positive integer such that $f^m(x)$ is a
periodic point. If $n$ is the period of $f^m(x)$ and $p \in
\big(n \mathbb{N} + l\big)^*$ for some $l < n$, then $f^p(x) =
f^l(f^{nj}(x))$ where $j=\min\{i:\, m\leq ni+l\}$.

\item[$(iii)$]  Suppose that the orbit of  $x$ is  infinite and $\omega_f(x)=\mathcal O_f(y)$
for some periodic point $y \in X$ with period $n$. If  $p,\,\,q\in(n\mathbb N+l)^*$ for some  $l<n$, then  $f^p(x)=f^q(x)$.

\item[$(iv)$]   $f^p(f^n(x))=f^n(f^p(x))$  for every  $n\in\nat$, $x\in X$ and every $p\in\nat^*$.
\end{itemize}
\end{lema}

The next statement is Lemma 2.4 of \cite{gru}.

\begin{lema}\label{aisladoperiodico2} Let $(X,f)$ be a dynamical system.  If
$\omega_f(x)$ is finite, then every point of $\omega_f(x)$ is
periodic. In particular, if $\omega_f(x)$ has an isolated point in
$\overline{\mathcal{O}_f(x)}$, then every point of $\omega_f(x)$
is periodic.
\end{lema}

The following lemma is a reformulation of Theorem 2.2  from  \cite{gru2016}.

\begin{lema} \label{Ellis finito} Let $(X,f)$ be a dynamical system, $x\in X$ and $m \in \mathbb{N}$.
Then  $|\mathcal O_f(x)|>m$  iff $f^m(x)\neq f^n(x)$ for every $m \neq n$.
In particular,  there is an integer  $M >0$ such that
$|\mathcal O_f(x)|< M$ for each $x\in X$ iff $E(X,f)$ is finite.
\end{lema}

We omit the proof of the following well-known result.

\begin{lema}\label{compactificacionorbitainfinita2} Let $(X,f)$ be a dynamical system such that $X$ has a dense subset consisting of isolated points.
If $X$ has a point with infinite orbit, then  $\{f^n:n\in \mathbb{N}\}$ is infinite and discrete in $E(X,f)$ and $f^n \neq f^p$ for every $(n,p) \in \mathbb{N} \times \mathbb{N}^*$. In addition, if the orbit of $w$ is dense in $X$, then  $X = \{ f^p(w) : p \in \beta(\mathbb{N}) \}$
 and thus  $X$ is a continuous image of $\beta(\nat)$.
\end{lema}

\section{Continuity of the $p$-iterates}

In the context of compact metric countable spaces,  we showed in  \cite[Th. 3.11]{gru} that if all points of $a\in X'$ are periodic,
then for $b\in X'$ either each $f^q$ is discontinuous at $b$ for every
$q\in\mathbb N^\ast$ or each $f^q$ is continuous at $b$ for all
$q\in\mathbb N^\ast$. An example where this assertion fails assuming that  all points of $X'$ are eventually periodic is also given  in \cite{gru}.  Here, our main task is to give
a necessary and sufficient condition on the space $X$ in order that all  $p$-iterates be discontinuous at a given point.
To have this done we shall need some preliminary lemmas.

\medskip

The next two results correspond to  Lemma 3.7 and  Theorem 3.8 of \cite{gru}, respectively.

\begin{lema}\label{orbita converge para X' periodicos}Let $(X,f)$ be a dynamical
system such that $X$ is a compact metric countable space and  every point of $X'$
is  periodic. If $x$ has an infinite orbit and  $y \in \omega_f(x)$  is fixed,
then $f^n(x)\longrightarrow y$.
\end{lema}

\begin{lema}\label{w_f(y)=orbita periodica} Let $(X,f)$ be a dynamical
system such that $X$ is a compact metric countable space and  every point of $X'$
is  periodic. For every   $x\in X$, there exists a periodic point $y \in X$ such that $\omega_f(x)=\mathcal O_f(y)$.
\end{lema}

The following two additional results are needed to establish our theorems.

\begin{lema}
\label{puntosquesaltan}
Let $f:X\rightarrow X$ be
a continuous function such that every accumulation point is fixed.
If $V$, $W$  are nonempty open sets such that $\overline V\cap \overline W=\emptyset$, then the set $\{x\in X: x\in V \text{ and  } f(x)\in W\}$
is finite.
\end{lema}

\proof Assume, towards a contradiction, that  $H= \{x\in X: x\in V \text{ y } f(x)\in W\}$ is infinite.
Since $X$ is compact and metric, there is a non constant  sequence $(a_n)_{n \in\mathbb N}$ in $H$ and $a\in \overline{V} \cap X'$
such that $a_{n}\rightarrow a$. As $f$ is continuous
and  $a$ is fixed, then   $f(a_{n})\rightarrow f(a) = a$ which implies that $a \in \overline{W}$.
But this is a contradiction. This shows that  $H$ is finite.
\endproof

\begin{lema}\label{puntoaislado} Let $(X,f)$ be a dynamical system such that $X$ is a compact metric countable space  and every accumulation point is periodic. Let  $p \in \mathbb{N}^*$ and $b \in X$ be an isolated point.  If there is  a sequence $(a_n)_{n \in \mathbb{N}}$ in $X$ such that $f^p(a_n) \to b$, then $b$ is periodic and $\mathcal{O}_f(b) = \omega_f(b) = \omega_f(a_n)$ for  every positive integer except  finitely many.
\end{lema}

\proof It follows directly from the assumption  that   $B=\{n\in\mathbb N:f^p(a_n)=b\}$ is cofinite and hence $b \in \omega_f(a_n) \subseteq \overline{\mathcal{O}_f(a_n)}$ for every $n\in B$.  By Lemma \ref{aisladoperiodico2}, we have that  $b$ is a  periodic point. So,
$\mathcal{O}_f(b) = \omega_f(b)$. By Lemma \ref{w_f(y)=orbita periodica},  for every $n \in B$ there is a periodic point $y_n \in X$ such that $\omega_f(a_n) = \mathcal{O}_f(y_n)$ and since $b \in \mathcal{O}_f(y_n)$, we conclude  that $\omega_f(b) = \mathcal{O}_f(b) = \mathcal{O}_f(y_n)  = \omega_f(a_n)$ for all $n \in B$.
\endproof

\begin{teo}\label{X' fija todas p-iteradas continuas} Let $(X,f)$ be a dynamical system such that $X$ is a compact metric countable space and
every accumulation point of $X$  is fixed. For every $a\in X'$, the following statements are equivalent:
\begin{enumerate}
  \item[$(1)$] There is $p\in\mathbb N^\ast$ such that $f^p$ is discontinuous at $a\in X'$.

  \item[$(2)$] There are a periodic point $b\in X\setminus \{a\}$  and a sequence $(a_n)_{n \in \mathbb{N}}$ in $X$ such that $a_n \to a$ and
  $b \in \overline{\mathcal{O}_f(a_n)}$ for all $n \in \mathbb{N}$.

  \item[$(3)$] $f^p$ is discontinuous at $a$, for each $p \in \mathbb N^\ast$.
\end{enumerate}
\end{teo}

\proof $(1)\Rightarrow (2)$.  Suppose that  $f^p$ is discontinuous  at $a\in X'$. Then,
          there is a nontrivial sequence $(a_n)_{n\in \mathbb N}$ in $X$ such that $a_n \to a$ and $f^p(a_n)$ does not converge to $a$.
          Since $X$ is compact and metric, there are a sequence of positive integers $(n_k)_{k \in \mathbb N}$
          and $b\in X \setminus \{a\}$ such that $f^p(a_{n_k}) \rightarrow b$.
          Assume, without loss of generality, that this subsequence  is
          $(a_n)_{n\in \mathbb N}$. In virtue of Lemma \ref{puntoaislado}, we only consider the case when  $b\in X'$. Since $a\neq b$,
          $f^p(a_n)\rightarrow b$ and $a_n\rightarrow a$,
          there are a clopen set  $V\in\mathcal{N}(a)$
          and $N \in\mathbb N$ such that  $b\notin V$, $a_n\in V$ y $f^p(a_n)\in X \setminus V$
          for each $n\geq N$.     For every  $n\geq N$ there is
          $d_n\in\mathcal{O}_f(a_n)\cap V$ such that $f(d_n)\in X \setminus V$, this is possible since  the set $A_n=\{m\in\mathbb N: f^m(a_n)\notin V\}\in p$ and hence it is infinite. Then,
          by Lemma \ref{puntosquesaltan}, we have that
          the set $B=\{d_n: n\geq N \}$ is finite. Hence, there
          exists $d\in B$ for which the set  $H=\{n\in\mathbb N: d=d_n\}$ is infinite.
          To finish the proof it suffices to show that $b\in \omega_f(a_n)$ for all $n \in H$.  Indeed, we consider two cases.
          Suppose first that $\mathcal{O}_f(d)$ is infinite, then
          there is $e\in X'$ such that $e\in\omega_f(d)$.
          As $e$ is fixed, by Lemma \ref{orbita converge para X' periodicos},
          $f^m(d)\rightarrow e$.  Analogously it is shown that $f^m(a_n)\xrightarrow[m \to \infty]{} e$ because of $e \in \omega_f(a_n)$, for each $n\in H$.
          Consequently, we obtain that $f^q(a_n)=e$ for each $q\in\mathbb N^\ast$ and for each  all $n\in H$. This implies that  $b=e$ since
          $f^p(a_n)\xrightarrow[n \in H]{} b$.  Thus $b\in \omega_f(a_n)$ for all $n\in H$.  For the second case, we assume that $\mathcal{O}_f(d)$ is finite. As $d \in \mathcal{O}_f(a_n)$ for each $n \in H$, we also have that the orbit of  $a_n$ is finite for all $n \in H$. By
          Lemma \ref{w_f(y)=orbita periodica}, we may choose   a periodic point  $e\in X$
          such that $\mathcal{O}_f(e)= \omega_f(d) = \{f^q(d):q\in \mathbb N^*\}=\{f^q(a_n):q\in \mathbb N^*\} = \omega_f(a_n)$
          for $n\in H$. Since    $\mathcal{O}_f(e)$ is finite
          and $f^p(a_n)\rightarrow b$, we obtain that $b\in \omega_f(a_n)$, for each $n\in H$.

\medskip

$(2) \Rightarrow (3)$.  We have to analyze two possible cases.

 Caso I. Suppose $b$ is isolated.  By assumption,  $b$ is periodic. Hence,  if $b \in \mathcal O_f(y)$, then
            $\mathcal{O}_f(y)$ is finite and so $\mathcal O_f(b) = \omega_f(y) = \{f^p(y): p\in \mathbb{N}^*\}$. Thus, from the hypothesis we obtain that
            $f^p(a_n) \in \omega_f(a_n) = \mathcal O_f(b)$ for every $n \in \mathbb{N}$ and  every $p \in \mathbb{N}^*$.  As $a$ is fixed and $b\neq a$, then  $a\nin \mathcal O_f(b)$ and the sequence $(f^p(a_n))_{n \in \mathbb{N}}$ cannot converge to $a$ for any $p\in\mathbb{ N}^\ast$. Thus $f$ is discontinuous at $a$.

            Caso II. Suppose $b\in X'$. Observe that
            if $b \in \overline{\mathcal O_f(y)}$ and $\mathcal O_f(y)$ is finite, then  $\omega_f(y)=\{b\}$.
            On the other hand, if $b \in \overline{\mathcal O_f(y)}$ and $\mathcal O_f(y)$ is infinite, then
            it follows from Lemma \ref{orbita converge para X' periodicos} that $\omega_f(y)=\{b\}$. Thus, we have that $f^p(a_n) = b$ for each $p\in\mathbb{ N}^\ast$ and for each $n \in \mathbb{N}$. Thus $f$ is discontinuous at $a$.

\medskip

$(3)\Rightarrow (1)$. It is evident.
\endproof

\begin{coro} Let $(X,f)$ be a dynamical system such that $X$ is a compact metric countable space and
every accumulation point of $X$  is fixed. If the orbit of every  isolated point of $X$ is finite, then every function of $E(X,f)$ is continuous.
\end{coro}

\proof Assume that there is $p\in \mathbb N^*$ such that $f^p$ is discontinuous at $a \in X'$. It follows from Theorem \ref{X' fija todas p-iteradas continuas}(2) that there are a periodic point $b\in X\setminus \{a\}$  and a sequence $(a_n)_{n \in \mathbb{N}}$ in $X$ such that $a_n \to a$ and
  $b \in \overline{\mathcal{O}_f(a_n)}$  for all $n \in \mathbb{N}$. Since each accumulation point of $X$ is fixed, we may assume that $a_n$ is isolated for every $n \in \mathbb{N}^*$, hence $b \in {\mathcal{O}_f(a_n)}$  for all $n \in \mathbb{N}$.  Now, let $V$ and $W$ be disjoint open sets such that $\overline{V} \cap \overline{W} = \emptyset$, $a \in V$ and $b \in W$.  Without loss of generality, we may assume that $a_n \in V$ for all $n \in \mathbb{N}$. Then we can find a sequence $(k_n)_{n \in \mathbb{N}}$ such that $f^{k_n}(a_{n}) \in V$ and $f^{k_n+1}(a_{n}) \in W$ for every $n \in \mathbb{N}$. We may assume that $f^{k_n}(a_{n}) \to c \in \overline{V}$, but this is impossible since $f^{k_n+1}(a_{n}) \to f(c) = c \in \overline{W}$. Therefore, $f^p$ is continuous.
  \endproof

The conclusion of Theorem  \ref{X' fija todas p-iteradas continuas} is not true if we replace the hypothesis  ``every  accumulation point is fixed'' by the hypothesis  ``every accumulation point is periodic''. Indeed, the next two examples witness
that the condition  $(2)$ of  Theorem \ref{X' fija todas p-iteradas continuas}
holds together  with either the assumption   ``$f^p$ is continuous for every $p\in \mathbb N^*$''   or  the assumption ``$f^p$ is discontinuous for all $p\in \mathbb N^*$''  .

\medskip

The phase space of the following dynamical systems  is the ordinal space $2\omega +1$ (two disjoint convergent sequences) which will be
identified with the following  subspace of $\mathbb{R}$:
$$
X=\{a_n:n\in\mathbb N\}\cup\{b_n:n\in\mathbb N\}\cup\{a,b\},
$$
where $a_n < a_{n+1} < a < b_n < b_{n+1}  < b$ for every $n \in \mathbb{N}$, $a_n \to a$ and $b_n \to b$.

\begin{ejem}\label{ejemplo1}  Define a function $f:X\rightarrow X$ as follows:
\begin{enumerate}
  \item[$a)$]   $f(a)=b$ and $f(b)=a$,

  \item[$b)$]   $f(a_n)=b_n$ for each $n\in\mathbb N$, and

  \item[$c)$]   $f(b_n)=a_{n+1}$ for each $n\in\mathbb N$.
\end{enumerate}
That is,
$$
a_{0}\rightarrow b_{0}\rightarrow a_{1}\rightarrow b_{1}\rightarrow a_{2}\rightarrow b_{2}
\rightarrow a_{3}\rightarrow b_{3}\cdots a_{n}\rightarrow b_{n}\rightarrow a_{n+1}\rightarrow b_{n+1}\rightarrow a_{n+2}\rightarrow b_{n+2}\rightarrow\cdots
$$
It is evident that $f$ is continuous. From the definition of $f$ we have, for all  $n, m \in\mathbb N$, the following:
\begin{enumerate}
  \item[$i)$]   $f^{2m}(a_n)=a_{n+m}$,

  \item[$ii)$]   $f^{2m+1}(a_n)=b_{n+m}$,

  \item[$iii)$]  $f^{2m}(b_n)=b_{n+m}$, and

  \item[$iv)$]  $f^{2m+1}(b_n)=a_{n+m+1}$.
\end{enumerate}
Conditions $i) - iv)$ imply the following:
\begin{enumerate}
\item[$(1)$] $f^p(a_n)=f^p(a)=a$ for every $p\in(2\mathbb N)^*$,

\item[$(2)$] $f^p(a_n)=f^p(b)=b$ for every $p\in(2\mathbb N+1)^*$,

\item[$(3)$] $f^p(b_n)=f^p(b)=b$ for every $p\in(2\mathbb N)^*$, and

\item[$(4)$] $f^p(b_n)=f^p(a)=a$ for every $p\in(2\mathbb N+1)^*$.
\end{enumerate}
Then we have that the accumulation points  $a$ and $b$ have period equal to $2$ and both  satisfy the second condition  of
Theorem \ref{X' fija todas p-iteradas continuas}. However, the function  $f^p$ is continuous for every $p\in\mathbb N^*$.
\end{ejem}

\begin{ejem}\label{ejemplo 2}  We define a function $f:X\rightarrow X$ as follows:
\begin{enumerate}
\item[$a)$]   $f(a)=b$ and $f(b)=a$,

\item[$b)$]   $f(a_0)=a_0$,

\item[$c)$]   $f(a_n)=b_{n-1}$ for every  $0 < n\in\mathbb N$, and

\item[$d)$]   $f(b_n)=a_{n}$, for every  $n\in\mathbb N$.
\end{enumerate}
That is,
$$
b_n\rightarrow a_n\rightarrow b_{n-1}\rightarrow a_{n-1}\rightarrow b_{n-2}\rightarrow a_{n-2}
\rightarrow b_{n-3}\rightarrow a_{n-3}\rightarrow\cdots b_{2}
\rightarrow a_{2}\rightarrow b_{1}\rightarrow a_{1}\rightarrow b_0\rightarrow a_0.
$$
It is not hard to prove that  $f$ is continuous. From the definition of $f$ we easily have  that for each $x\in X\setminus\{a,b\}$  there exists $n\in\mathbb N$
such that $f^n(x)=a_0$. Hence, $f^p(x)=a_0$ for every $x\in X\setminus\{a,b\}$ and every
$p\in\mathbb N^*$. Since $f^p(a)=a$ and $f^p(b)=b$ for all $p\in (2\mathbb N)^*$; and
$f^p(a)=b$ and $f^p(b)=a$ for all $p\in (2\mathbb N+1)^*$, we conclude that $f^p$ is
discontinuous at $a$ and $b$ for each $p\in\mathbb N^*$. Finally, observe that
$a_0\in \overline{\mathcal{O}_f(a_n)}$ for all $n \in \mathbb{N}$, so $(2)$ of
Theorem \ref{X' fija todas p-iteradas continuas} is satisfied.
\end{ejem}

\medskip

In the next theorem, we show a generalization of $(2) \Rightarrow (3)$ in  Theorem  \ref{X' fija todas p-iteradas continuas}.

\begin{teo}
\label{nuevo2}
 Let $(X,f)$ be a dynamical system such that $X$ is a compact
metric countable space and every accumulation point of $X$  is periodic.
Let $a\in X'$. If there are a sequence
$(a_n)_{n \in \mathbb{N}}$ in $X$ and a periodic point $b\in X\setminus \mathcal O_f(a)$ such that $a_n \to a$ and
$b \in \overline{\mathcal{O}_f(a_n)}$ for every $n \in \mathbb{N}$, then
$f^p$ is discontinuous at $a$ for each $p\in \mathbb N^\ast$.
\end{teo}

\proof   First  we show a particular case. Suppose that  $b \in \mathcal O_f(a_n)$ for infinitely many $n \in \mathbb{N}$.   Hence,   from the periodicity of $b$, we have that     $\mathcal O_f(b) = \omega_f(b)= \omega_f(a_n) = \{f^p(a_n): p\in \mathbb{N}^*\}$ is finite for infinitely many $n \in \mathbb{N}$. Since $\mathcal{ O}_f(a) \cap \mathcal {O}_f(b) = \emptyset$, the  sequence $(f^p(a_n))_{n \in \mathbb{N}}$ cannot converge to $f^p(a) \in \mathcal{O}_f(a)$ for any $p\in\mathbb{ N}^\ast$. Thus $f^p$ is not continuous at $a$ for any  $p\in\mathbb{ N}^\ast$.

For the proof of the theorem we consider two cases: (i) Suppose  $b$ is isolated. Then $b \in  \mathcal O_f(a_n)$ for every $n \in \mathbb{N}$ and we are done from the result we proved above.

(ii) Assume that $b\in X' \setminus \big(\bigcup_{n \in \mathbb{N}}\mathcal{O}_f(a_n)\big)$. From the result proved above, we can assume  that $\mathcal O_f(a_n)$ is infinite for every  $n \in \mathbb{N}$. Then, we must have that $b \in \omega_f(a_n)$ for every $n \in \mathbb{N}$. Now,  in virtue of Lemma \ref{w_f(y)=orbita periodica}, for every $n \in \mathbb{N}$ there is a periodic point  $y_n \in X$ such that $\omega_f(a_n)=\mathcal{O}_f(y_n)$. Since $b$ is periodic and $b \in  \mathcal{O}_f(y_n)$, then $\mathcal{O}_f(b) = \mathcal{O}_f(y_n)$ for all $n \in \mathbb{N}$.  Thus $f^p(a_n)\in \mathcal{O}_f(b)$ for all $n$ and all $p\in\mathbb{ N}^\ast$. Since $\mathcal{ O}_f(a) \cap \mathcal {O}_f(b) = \emptyset$, we conclude as before that   $f^p$ is not continuous at $a$ for any  $p\in\mathbb{ N}^\ast$.
\endproof

In the next result we slightly modify the proof of  the previous theorem to get an interesting statement.

\begin{teo} Let $(X,f)$ be a dynamical system such that $X$ is a compact metric countable space and
every accumulation point of $X$  is periodic. Let $a \in X'$ and $(a_n)_{n \in \mathbb{N}}$ be  a sequence in $X$ such that
$f^p(a_n) \to f^p(a)$, for some $p \in \mathbb{N}^*$. Suppose  $b \in X$ is a periodic point  and $b \in \bigcap_{n \in \mathbb{N}} \overline{\mathcal{O}_f(a_n) }$, then $b \in\mathcal{O}_f(a)$.
\end{teo}

\proof  Let $a\in X'$, $(a_n)_n$  in $X$ and $p\in \mathbb{N}^*$ as in the hypothesis. First, suppose that $b$ is isolated. Then, $b \in \bigcap_{n \in \mathbb{N}}\mathcal{O}_f(a_n)$ and so $\mathcal{O}_f(b) = \omega_f(b) = \omega_f(a_n)$ for every $n \in \mathbb{N}$. Since $f^p(a_n) \to f^p(a) \in \mathcal{O}_f(a)$ and $f^p(a_n) \in  \omega_f(a_n) = \mathcal{O}_f(b)$ for each $n \in \mathbb{N}$, we must have that $\mathcal{O}_f(a) \cap \mathcal{O}_f(b) \neq \emptyset$ and so $b \in \mathcal{O}_f(a)$. Now, suppose that $b \in X'$. Notice that if  $b \in \mathcal{O}_f(a_n)$ for some $n \in \mathbb{N}$, then $\mathcal{O}_f(b) = \omega_f(b) =  \omega_f(a_n)$. Thus we have that  $b \in \omega_f(a_n)$ for all $n \in \mathbb{N}$.
By Lemma \ref{w_f(y)=orbita periodica}, there exists a periodic point $y_n \in X$ such that  $\omega_f(a_n)=\mathcal O_f(y_n)$. Since $b \in \omega_f(a_n)$, then $\omega_f(a_n) = \mathcal{O}_f(b)$ for all $n\in \mathbb N$. Thus, we have shown that $f^p(a_n) \in \mathcal{O}_f(b)$ for every $n \in \mathbb{N}$. Then, as before, we have $\mathcal{O}_f(a) \cap \mathcal{O}_f(b) \neq \emptyset$  and thus  $b \in \mathcal{O}_f(a)$.
\endproof

Concerning the last theorem we have the following example.

\begin{ejem} We consider again the  countable ordinal space $2\omega + 1$ identified with the subspace $X$ of $\mathbb{R}$ from above.
Define the function $f:X\rightarrow X$ as follows:
\begin{enumerate}
  \item[$a)$]   $f(a)=b$ and $f(b)=a$,

  \item[$b)$]   $f(a_0)=b$,

  \item[$c)$]   $f(a_n)=b_{n-1}$ for each $n>0$, and

  \item[$d)$]   $f(b_n)=a_{n}$  for each $n\in\mathbb N$.
\end{enumerate}
The function $f$ is evidently continuous. From the definition we can see that
$$
b_n\rightarrow a_n\rightarrow b_{n-1}\rightarrow a_{n-1}\rightarrow b_{n-2}\rightarrow a_{n-2}
\rightarrow\cdots b_{2}\rightarrow a_{2}\rightarrow b_{1}\rightarrow a_{1}\rightarrow b_0\rightarrow a_0\rightarrow b.
$$
Hence, all points are eventually periodic. Besides, we have the following properties:
\begin{enumerate}
\item[$i)$]  For every $x\in X$ there is  $n \in\mathbb{N}$ such that  $f^n(x)=b$.

\item[$ii)$]   $f^p(x)=b$ for each $x\in X\setminus\{a,b\}$ and for each $p\in\mathbb{N}^*$.

\item[$iii)$]   $f^p(a)=a$  for all $p\in (2\mathbb N)^*$.

\item[$iv)$]   $f^p(a)=b$  for all $p\in (2\mathbb N +1)^*$.

\item[$v)$]   $f^p(b)=a$  for all $p\in (2\mathbb N +1)^*$.
\end{enumerate}
Thus,  we have that $f^p$ is discontinuous at $a$ for all $p\in (2\mathbb N)^*$ and we also have that $b \in \mathcal{O}_f(x)$ for every $x \in X$. Moreover, $f^p$ is discontinuous at $b$ for all $p\in (2\mathbb{N} + 1)^*$
\end{ejem}

Theorem \ref{X' fija todas p-iteradas continuas}  suggests the following problem.

\begin{problema} Let $(X,f)$ be a dynamical system such that $X$ is a compact metric countable space such that
 every accumulation point of $X$  is periodic. Find a necessary and sufficient condition on a point $a \in X$, like in Theorem \ref{X' fija todas p-iteradas continuas}, in order that  $f^p$ is discontinuous at $a$ for each $p \in \mathbb N^\ast$.
\end{problema}

\section{Transitive dynamical systems}

Continuing the work presented in \cite{gru2016}, in this section we focus our attention on transitive dynamical systems. We recall some questions from that paper.  The first one is whether $E(X,f)$ is countable for every transitive system $(X,f)$ (see Questions 4.6 in \cite{gru2016}). For instance, this happens when every function in $E(X,f)$ is continuous (see \cite[Theorem 2.9, Theorem 3.3]{gru2016}). Below we extend this result.
The second question is  whether or not   is there a continuous function $f: X \to X$ such that $E(X,f)$ is homeomorphic to $Y$, where $X$ and $Y$ are arbitrary compact metric countable  spaces? (see Questions 4.8 in \cite{gru2016}).

We need the following result from \cite{gru2016} (see Lemma 3.1 and Theorem 2.3).

\begin{lema}
\label{p(f)finitoasterisco}
Let $(X,f)$ be a dynamical system.

(i) $E(X,f)\setminus\{f^n:\;n\in \mathbb N\}$ is finite
iff there is $m \in \mathbb{N}$ such that $|\omega_f(x)|\leq m$ for each $x\in X$.

(ii) If  $(X,f)$ is transitive and $y\in X'$, then $f(y)\in X'$.

(iii)  If $w$ has a dense orbit, then $w$  is isolated.
\end{lema}

\begin{teo}
\label{homeo}
Let $(X,f)$ be a transitive dynamical system where $X$ is a compact metri\-zable countable space. If either $f^p$ is continuous for each $p\in \mathbb{N}^*$, or   $|\omega_f(x)|=1$ for each $x \in X'$, then  $E(X,f)$ is homeomorphic to $X$.
\end{teo}

\proof Let $w$ be a point of $X$ whose  orbit is dense in $X$. According to Lemma \ref{p(f)finitoasterisco} the point $w$ must be isolated. Now,  consider the function $h:E(X,f)\rightarrow X$ defined by
$h(f^p)=f^p(w)$, for every $p \in \beta(\mathbb{N})$.  This function is continuous since it is the restriction of the projection map $\pi_w : X^X \to X$ to $E(X,f)$.  It follows from  Lemma \ref{compactificacionorbitainfinita2} that $h$ is surjective since $X= \{ f^p(w) : p \in \beta(\mathbb{N}) \}$.  To prove the theorem  it  suffices to show that the function
$h$ is injective. To have this done, first observe that $f^p(f^n(w))=f^n(f^p(w))$, for all $n \in \mathbb{N}$ and for all  $p \in \beta(\mathbb{N})$, and the orbit of $w$ is the collection of isolated points of $X$. Hence, we obtain that $f^p(x) = f^q(x)$ for every isolated point $x \in X$ whenever $f^p(w) = f^q(w)$ for some $p, q \in \beta(\mathbb{N})$.

\medskip

Suppose first that $f^p$ is continuous, for each $p\in \mathbb{N}^*$. Then, if $p, q\in \mathbb{N}^*$ and  $f^p$ and $f^q$ agree on all the  points of a dense orbit, then we obtain that $f^p = f^q$. This shows $h$ is injective.

\medskip

Now, assume that  $|\omega_f(x)|=1$ for each $x \in X'$. Then for every  $x \in  X'$ there is $z_x \in X'$ such that $\omega_f(x) = \{z_x\}$ and hence we obtain that  $f^p(x) = z_x$ for all $p \in \mathbb{N}^*$. Thus we have that $f^p=f^q$ iff $f^p(w)=f^q(w)$ for all $p,q\in \beta(\nat)$. So, $h$ is injective.

\medskip

 Therefore, in both cases, $h$ is a homeomorphism between $E(X,f)$ and $X$.
\endproof

Question \ref{final} is related to the second condition of the previous theorem.

\begin{coro} Let $(X,f)$ be a transitive dynamical system  where $X$ is a compact metri\-zable countable space.
If $f^p$ is continuous for each $p\in \mathbb{N}^*$, then the set $P_f$ is finite.
\end{coro}

\proof By Theorem \ref{homeo}, we know that $E(X,f)$ is countable. If $P_f$ were infinite, by  Theorem 2.7 of \cite{gru2016}, then  $E(X,f)$ would be homeomorphic to the Cantor set $2^\mathbb{N}$. So, $P_f$ is finite.
\endproof

As we have seen in the proof of the  previous corollary, if $X$ is a compact metrizable space and  $E(X,f)$ is countable, then $P_f$ must be finite. Next we shall estimate the cardinality of the Ellis semigroup under some restriction on the $\omega$-sets.

\begin{teo}\label{countable} Let $(X,f)$ be a transitive dynamical system such that
$X$ is a compact metric coun\-table  space. If there is   $m \in \mathbb{N}$ such that  $|\omega_f(x)|\leq m $ for every  $x\in X'$, then
$E(X,f)$ is countable.
\end{teo}

\proof Let $E^*=E(X,f)\setminus\{f^n:\;n\in \mathbb N\}$. It suffices to show that $E^*$ is countable.
First notice that $E^*$ is equal to $\{f^p:\; p\in \mathbb{N}^*\}$ by Lemma \ref{compactificacionorbitainfinita2}. Since $X'$ is $f$-invariant (by Lemma \ref{p(f)finitoasterisco}(ii)), then  $(X',f\restriction X')$ is a well defined dynamical system. By Lemma \ref{p(f)finitoasterisco}(i), $E(X',f\restriction X')$  is finite.
 Let  $w \in X$ with a dense orbit.   Consider the function  $\varphi: E^*\to E(X', f\restriction X')\times X'$ given by
$\varphi (g)=(g\restriction X', g(w))$. It suffices to show that  $\varphi$ is injective.  In fact, let $g\in E^*$, then $g=f^p$ for some $p\in \mathbb{N}^*$. Notice that  every isolated point is of the form $f^l(w)$ for some $l\in \mathbb N$. Thus $g(f^l(w))=f^p(f^l(w))= f^l(f^p(w))=f^l(g(w))$.
Therefore $g$ is completely determined by $g\restriction X'$ and $g(w)$.  Hence $\varphi$ is injective.
\endproof

Our next example satisfies the second conditions of Theorem \ref{homeo} and  all the  $p$-iterates are discontinuous, for $p \in \mathbb{N}^*$.

\medskip

The phase space of the next two examples is going to be   $\omega^3+1$ which, for our convenience,  will be identified with the following subspace of $\mathbb{R}$:
$$
X= \{d_{i,j,k}:i, j, k \in\mathbb N\}
\cup \{d_{j,k} :j, k \in \mathbb N \}\cup\{d_k:k\in\mathbb N\}\cup \{d\},
$$
where $(d_k)_{k \in \mathbb N}$ is a strictly increasing sequence such that $d_{k} \xrightarrow[k \to \infty]{} d$; $(d_{j,0})_{j \in \mathbb N}$ is a strictly increasing sequence  contained in $(-\infty,d_0)$ such that  $d_{j,0} \xrightarrow[j \to \infty]{} d_0$; for each positive $k \in \mathbb{N}$, the sequence
$(d_{j,k})_{j \in \mathbb N}$ is strictly increasing,  it is contained in $(d_{k-1},d_k)$ and  $d_{j,k} \xrightarrow[j \to \infty]{} d_k$; $D_{0,0}:=\{d_{i,0,0}: i\in\mathbb N\}$ is a strictly increasing sequence such that  $d_{i,0,0} \xrightarrow[i \to \infty]{} d_{0,0}$ and it is contained in $(-\infty,d_{0,0})$; $D_{0,k} =\{d_{i,0,k}:i\in\mathbb N\}$ is a strictly increasing sequence with $d_{i,0,k} \xrightarrow[i \to \infty]{} d_{0,k}$ and contained in $(d_{k-1},d_{0,k})$  for each $k \in \mathbb N \setminus \{0\}$ $D_{j,k}:=\{d_{i,j,k}:i \in\mathbb N\}$ is a strictly increasing sequence with   $d_{i,j,k} \xrightarrow[i \to \infty]{} d_{j,k}$ and contained in $(d_{j-1,k},d_{j,k})$  for every $j \in \mathbb N \setminus \{0\}$ and for every $k \in \mathbb{N}$.

\medskip

We are ready to describe our first  example.

\begin{ejem} There is a continuous function $f:\omega^3+1 \rightarrow \omega^3+1$ such that:
\begin{enumerate}
\item[$(1)$] $\mathcal{O}_f(d_{0,0,0})$ is dense.

\item[$(2)$]  The points $d_{c_n -1,0}$ and $d_{c_n +1,0}$ have infinite  orbits, where  $c_n=2+3n$ for every $n \in \mathbb{N}$.

\item[$(3)$]  $|\omega_f(x)| = 1$ for every $x \in (\omega^3+1)'$.

\item[$(4)$]  $E(\omega^3+1,f)$ is homeomorphic to $\omega^3+1$.

\item[$(5)$]  $f^p$ is discontinuous for all $p\in\mathbb N^*$.
\end{enumerate}
Our  function $f$ is defined as follows:
\begin{enumerate}
\item[$a$)]  $f(d)=d$ and  $f(d_{n})=d_{n}$ for each $n\in\mathbb N$.

  \item[$b)$]  $f(d_{0,0})=d$, $f(d_{1,0})=d_{3,0}$, $f(d_{2,0})=d_{0,0}$ and  $f(d_{0,1})=d_0$.

  \item[$c)$]  $f(d_{c_n-1,0})=d_{c_{n-1}-1,0}$ for each $n>0$.

  \item[$d)$]  $f(d_{c_n,0})=d_{c_{n-1},0}$ for each $n>0$.

  \item[$e)$]  $f(d_{c_n+1,0})=d_{c_{n+1}+1,0}$ for each $n\in\mathbb N$.

  \item[$f)$]  $f(d_{j,k})=d_{j-1,k}$ for each $j>0$ and $k>0$.

  \item[$g)$]  $f(d_{0,k})=d_{k-1}$ for each $k>0$.

  \item[$h)$]  $f(d_{i,j,k})=d_{i+1,j-1,k}$ for each $i\in\mathbb N$ and $j,k>0$.

  \item[$i)$]  $f(d_{i,0,2})=d_{0,i,1}$ for each $i\in\mathbb{N}$.

  \item[$j)$]  $f(d_{i,0,k})=d_{0,i+1,k-1}$ for each $i\in\mathbb N$ and $k>2$.

  \item[$k)$]  $f(d_{i,0,1})=d_{0,c_{i+1}-1,0}$ for each $i\in\mathbb{N}$.

  \item[$l)$]  $f(d_{i,c_0,0})=d_{i+1,0,0}$ for each $i\in\mathbb N$.

  \item[$m)$]  $f(d_{i,c_n,0})=d_{i+1,c_{n-1},0}$ for each $i\in\mathbb N$ and $n>0$.

  \item[$n)$]  $f(d_{i,c_n-1,0})=d_{i+1,c_{n-1}-1,0}$ for each $i\in\mathbb N$ and $n>0$.

  \item[$o)$]  $f(d_0, c_{n+1},0)=d_{0,c_n,0}$ for each $n\in\mathbb N$.

  \item[$p)$]  $f(d_{i,c_n+1,0})=d_{i-1,c_{n+1}+1,0}$ for each $i>0$ and $n\in\mathbb N$.

  \item[$q)$]  $f(d_{i,4,0})=d_{i,1,0}$ for each $i\in\mathbb N$.

  \item[$r)$]  $f(d_{i,1,0})=d_{i,3,0}$ for each $i\in\mathbb N$.

  \item[$s)$]  $f(d_{i,0,0})=d_{0,0,i+2}$ for each $i\in\mathbb N$.

  \item[$t)$]  $f(d_{i,2,0})=d_{i+1,0,0}$ for each $i\in\mathbb N$.
\end{enumerate}
In the next diagram, we can see the behavior of the orbits on the isolated points.
$$
d_{0,0,0}\rightarrow d_{0,0,2}\rightarrow d_{0,0,1}\rightarrow d_{0,4,0}\rightarrow d_{0,1,0}
\rightarrow d_{0,3,0}\rightarrow d_{0,2,0}\rightarrow d_{1,0,0}
$$
\smallskip
$$
d_{1,0,0}\rightarrow d_{0,0,3}\rightarrow d_{0,1,2}\rightarrow d_{1,0,2}
\rightarrow d_{0,1,1}\rightarrow d_{1,0,1}\rightarrow d_{0,7,0}\rightarrow d_{1,4,0}
$$
$$
\rightarrow d_{1,1,0}\rightarrow d_{1,3,0}\rightarrow d_{0,6,0}\rightarrow d_{0,5,0}
\rightarrow d_{1,2,0}\rightarrow d_{2,0,0}
$$
\smallskip
$$
d_{2,0,0}\rightarrow d_{0,0,4}\rightarrow d_{0,1,3}\rightarrow d_{1,0,3}\rightarrow d_{0,2,2}
\rightarrow d_{1,1,2}\rightarrow d_{2,0,2}\rightarrow d_{0,2,1}\rightarrow d_{1,1,1} \rightarrow d_{2,0,1}\rightarrow d_{0,10,0}
$$
$$
\rightarrow d_{1,7,0}\rightarrow d_{2,4,0}\rightarrow d_{2,1,0}\\
\rightarrow d_{2,3,0}\rightarrow d_{1,6,0}\rightarrow d_{0,9,0}
\rightarrow d_{0,8,0}\rightarrow d_{1,5,0}
\rightarrow d_{2,2,0}\rightarrow d_{3,0,0}
$$
\smallskip

$$
d_{3,0,0}\rightarrow d_{0,0,5}\rightarrow d_{0,1,4}\rightarrow d_{1,0,4}\rightarrow d_{0,2,3}
\rightarrow d_{1,1,3}\rightarrow d_{2,0,3}\rightarrow d_{0,3,2}\rightarrow d_{1,2,2}\rightarrow
d_{2,1,2}\rightarrow d_{3,0,2}\rightarrow d_{0,3,1}
$$
$$
\rightarrow d_{0,2,2} \rightarrow d_{1,1,4} \rightarrow d_{2,1,3} \rightarrow d_{1,1,2}\rightarrow d_{1,0,4}\rightarrow d_{2,0,3} \rightarrow d_{2,0,2}\rightarrow d_{0,3,1}\rightarrow d_{1,2,1}
$$
$$
\rightarrow d_{2,1,1}\rightarrow d_{3,0,1}\rightarrow d_{4,0,0}
$$
$$\vdots \ \ \ \ \ \ \ \ \ \ \ \ \vdots \ \ \ \ \ \ \ \ \ \ \ \ \vdots$$
$$
d_{n,0,0}\rightarrow d_{0,0,n+2}\rightarrow d_{0,1,n+1}\cdots d_{1,0,n+1}\cdots\rightarrow d_{n,0,1}\rightarrow d_{0,c_{n+1}-1,0}\rightarrow\cdots
$$
$$
\rightarrow d_{n,2,0}\rightarrow d_{n+1,0,0}\cdots
$$
$$
\vdots \ \ \ \ \ \ \ \ \ \ \ \ \vdots \ \ \ \ \ \ \ \ \ \ \ \ \vdots
$$
\noindent
The following properties are satisfied.

\medskip

\begin{enumerate}
 \item[$i)$]     $f(D_{j,k})=D_{j-1,k}$ for each $j>0$ and  $k>0$.
 \item[$ii)$]   $f(D_{c_0,0})=D_{0,0}\setminus\{d_{0,0,0}\}$ and
 $f(D_{c_n,0})=D_{c_{n-1},0}\setminus\{d_{0,c_{n-1},0}\}$ for each $n>0$.
 \item[$iii)$]   $f(D_{c_n+1,0})= D_{c_{n+1}+1,0}\cup\{d_{0,c_n,0}\}$ for each $n\in\mathbb N$.
 \item[$iv)$]   $f(D_{c_n-1,0})=D_{c_{n-1}-1,0}\setminus\{d_{0,c_{n-1}-1,0}\}$ for each $n\neq0$.
 \item[$v)$]    $f(D_{0,1})=\{d_{0,c_{n+1}-1,0}:n\in\mathbb N\}$.
 \item[$vi)$]   $f(D_{0,0})=\{d_{0,0,k+2} : k\in\mathbb N\}$.
 \item[$vii)$]  $f(D_{1,0})=D_{3,0}$ and $f(D_{2,0})=D_{0,0}$.
 \item[$viii)$] $f(D_{0,2})=\{d_{0,j,1}:j\in\mathbb N\}$ and
 $f(D_{0,k})=\{d_{0,j+1,k-1}:j\in\mathbb N\}$ for each $k>2$.
\end{enumerate}

\medskip

Observe that clauses $a)$ and $i)$ imply that $f$ is continuous at $d$, at
$d_i$  for each $i>0$ and at $d_{ij}$ for every $i,j>0$. For $j>1$,  $f$ is
continuous at $d_{0j}$ by clauses $g)$ and $viii)$. By
$d)$, $e)$, $ii)$, $iii)$ and $iv)$, we have that  $f$ is continuous at $d_0$ and at
$d_{i,0}$ for every $i>2$. By  $b)$, $v)$, $vi)$ and $vii)$,  $f$ is
continuous at the points $d_{0,0}$,  $d_{10}$,  $d_{20}$ and $d_{0,1}$. Therefore, $f$ is continuous.

It is evident that  the orbit $\mathcal{O}_f(d_{0,0,0})$ is dense. Also it is evident that the points $d_{c_n -1,0}$ and $d_{c_n +1,0}$ have infinite  orbits, for every $n \in \mathbb{N}$.
The following relationships  follows directly from the definition:
\begin{enumerate}
  \item[$I)$]  $f^p(d)=d$ and  $f^p(d_{n})=d_{n}$ for each $n\in\mathbb N$ and $p \in \mathbb{N}^*$.
 \item[$II)$]  $f^p(d_{0,0})=d$,  $f^p(d_{2,0})=d$ and  $f^p(d_{0,1})=d_0$.
 \item[$III)$]  $f^p(d_{c_n,0})=d$ for each $n>0$.
  \item[$IV)$]  $f^p(d_{c_n+1,0}) =  d_0 = f^p(d_{c_{n}-1,0})$ for each $n\in\mathbb N$.
    \item[$V)$]  $f^p(d_{j,k})=d_{k-1}$ for each $j>0$ and $k>0$.
 \end{enumerate}
 All these properties imply clause $(3)$, condition $(4)$ follows from Theorem \ref{homeo} and the last condition $(5)$ follows from clauses $III)$ and $IV)$.
 \end{ejem}

In the next example, we show that there exists a continuous function $f$ such that the dynamical system $(\omega^3+1,f)$ is transitive  and has a sequence  of accumulation points  with arbitrarily large period. This example differs from Example 4.1 of \cite{gru} since this dynamical system is transitive and in the other one all points have finite orbit. Notice also that,  in the example below, all $p$-iterates are discontinuous.

\begin{ejem} There is a function $f:\omega^3+1 \rightarrow \omega^3+1$ such that
\begin{enumerate}
\item[$(1)$] $\mathcal{O}_f(d_{0,0,0})$ is dense.

\item[$(2)$]  All points of $(\omega^3+1)'$ have finite orbit.

\item[$(3)$]  The accumulation point $d_{a_n}$ has period $n+2$, where  $a_0=0$ and $a_n=a_{n-1}+ n+1$, for every $n\in \mathbb N$.

\item[$(4)$]  $E(\omega^3+1,f)$ is homeomorphic to $2^\mathbb N$.

\item[$(5)$]   $f^p$ is discontinuous for all  $p\in\mathbb N^*$.
\end{enumerate}
 The function  $f$ is defined as follows:
\begin{enumerate}
\item[$a)$] $f(d)=d$.

  \item[$b)$] $f(d_{a_n})=d_{a_n+n+1}= d_{a_{n+1}-1}$ for each $n\in\mathbb N$.

  \item[$c)$] $f(d_{a_n+k})=d_{a_n+k-1}$ for each $n\in \mathbb N$ and $0<k\leq n+1$.

  \item[$d)$] $f(d_{0,0})=d$ and $f(d_{0,a_n})=d_{a_{n-1}}$ for each $n>0$.

  \item[$e)$] $f(d_{i,a_n})=d_{i-1,a_n + n+1} = d_{i-1,a_{n+1} -1}$ for each $n\in\mathbb N$ and $i>0$.

  \item[$f)$] $f(d_{i,a_n+k})=d_{i,a_n+k-1}$ for each $i\in\mathbb N$, $n\in \mathbb N$ and $0<k\leq n+1$.

  \item[$g)$] $f(d_{0,0,0})=d_{0,0,1}$ and $f(d_{n,0,0})=d_{0,0,a_n}$  for each $n>0$.

  \item[$h)$] $f(d_{i,j,0})=d_{i+1,j-1,1}$ for each $i\in\mathbb N$ and $j>0$.

  \item[$i)$] $f(d_{i,0,1})=d_{i+1,0,1}$ for each $i\in\mathbb N$.

  \item[$j)$] $f(d_{i,j,1})=d_{i,j,0}$ for each $i\in\mathbb N$ and $j>0$.

  \item[$k)$] $f(d_{i,j,a_n})=d_{i,j-1,a_n+n+1}$ for each $i\in\mathbb N$ and $j,n>0$.

  \item[$l)$] $f(d_{i,0,a_n})=d_{0,i+1,a_n-1}$ for each $i\in\mathbb N$ and $n>0$.

  \item[$m)$] $f(d_{i,j,k})=d_{i,j,k-1}$ for each $i,j\in\mathbb N$ and $k\notin H$
            and $k\notin \{a_n+1:n\in\mathbb N\}$.

  \item[$n)$] $f(d_{i,j,k})=d_{i,j,k-1}$ for each $i\in\mathbb N$, $j>0$ and $k\in \{a_n+1:n>0\}$.

  \item[$o)$] $f(d_{i,0,k})=d_{i+1,0,k-1}$ for each $i\in\mathbb N$ and $k\in \{a_n+1:n>0\}$.
\end{enumerate}
To have some idea about the orbits  we describe some of them in the next diagram:

\medskip

\noindent $\underline{d_{0,0,0}} \rightarrow d_{0,0,1}\rightarrow \underline{d_{1,0,0}}
\rightarrow d_{0,0,2}\rightarrow d_{0,1,1}\rightarrow d_{0,1,0}
\rightarrow d_{1,0,1}\rightarrow \underline{d_{2,0,0}} \rightarrow d_{0,0,5}\rightarrow d_{0,1,4}\rightarrow d_{0,1,3}\rightarrow d_{0,1,2}
\rightarrow d_{0,0,4}\rightarrow d_{0,0,3}\rightarrow d_{0,0,2}\rightarrow d_{1,0,2}
\rightarrow d_{0,2,1}\rightarrow d_{0,2,0}\rightarrow d_{1,1,1}\rightarrow d_{1,1,0}\rightarrow d_{2,0,1}
\rightarrow \underline{d_{3,0,0}} \rightarrow d_{0,0,9}\rightarrow d_{0,1,8}\rightarrow d_{0,1,7}\rightarrow d_{0,1,6}
\rightarrow d_{0,1,5}\rightarrow d_{0,0,8}\rightarrow d_{0,1,7}\rightarrow d_{0,1,6}\rightarrow
d_{1,0,5}\rightarrow d_{0,2,4}\rightarrow d_{0,2,3}\rightarrow d_{0,2,2}\rightarrow d_{1,1,4}
\rightarrow d_{2,1,3}\rightarrow d_{1,1,2}\rightarrow d_{1,0,4}\rightarrow d_{2,0,3}
\rightarrow d_{2,0,2}\rightarrow d_{0,3,1}\rightarrow d_{0,3,0}\rightarrow d_{1,2,1}
\rightarrow d_{1,1,0}\rightarrow d_{3,1,1}\rightarrow d_{3,1,0}\rightarrow d_{3,0,1}\rightarrow d_{4,0,0} \rightarrow \underline{d_{0,0,14}} \rightarrow \ldots$
$$\vdots \ \ \ \ \ \ \ \ \ \ \ \ \vdots \ \ \ \ \ \ \ \ \ \ \ \ \vdots$$
\noindent $\rightarrow  \underline{d_{n,0,0}} \rightarrow d_{0,0,a_n}\rightarrow d_{0,1,a_{n-1}}\cdots d_{0,1,a_{n-1}}\rightarrow d_{0,0,a_n-1}\cdots\rightarrow d_{1,0,a_{n-1}}\rightarrow d_{0,2,a_{n-1}-1}\rightarrow\cdots
\rightarrow d_{0,2,a_{n-2}}\rightarrow d_{n,0,1}\rightarrow \underline{d_{n+1,0,0}} \rightarrow  \cdots$\\
$$\vdots \ \ \ \ \ \ \ \ \ \ \ \ \vdots \ \ \ \ \ \ \ \ \ \ \ \ \vdots$$
\noindent Notice  that:
\begin{enumerate}
  \item[i.]    $f(D_{0,0})=\{d_{0,0,a_n} : n\in\mathbb N \}$.
  \item[ii.]   $f(D_{j,0})=D_{j-1,1}\setminus\{d_{0,j-1,1}\}$ for each $j>0$.
  \item[ii.]   $f(D_{0,1})=D_{0,0}\setminus\{d_{0,0,0}\}$ and  $f(D_{ij})= D_{(i-1)j}$ for each $i, j>0$.
  \item[iii.]  $f(D_{j,1})= D_{j,0}$ for each $j>0$.
  \item[iv.]   $f(D_{0,a_n})= \{d_{0,j+1,a_n-1}:j\in\mathbb N\}$ for each $n>0$.
  \item[v.]    $f(D_{j,a_n})=D_{j-1,a_{n}+n+1}\setminus \{d_{0,j-1,a_n+n+1}\}$ for each $j$, $n>0$.
  \item[vi.]   $f(D_{j,k})=D_{j,k-1}$ for each $j\in\mathbb N$ and  $k\notin H \cup \{a_n+1:n\in\mathbb N\}$.
  \item[vii.]  $f(D_{j,k})=D_{j,k-1}$ for each $j>0$ and $k\in\{a_n+1:n\in\mathbb N\}$.
  \item[viii.] $f(D_{0,k})=D_{0,k-1}\setminus\{d_{0,0,k-1}\}$ for each $k\in\{a_n+1:n\in\mathbb N\}$.
\end{enumerate}

\bigskip

By conditions $a)$,  $b)$, $c)$, $iv)$, $v)$, $vi)$, $vii)$ and  $viii)$, we have that $f$ is continuous at $d$, $d_k$ for every $k>0$;
 $f$ is continuous at $d_0$ by conditions $b)$ and $ii)$;  and $f$ is continuous at $d_{ij}$,by conditions $d)$, $e)$, $f)$ and identities from $i)$ to $viii)$ for each
 $i, j\in\mathbb N$.  Consider the sequence  $(d_{i0})_{i\in\mathbb N}$  that converges to $d_0$. By conditions $d)$, $e)$ and $f)$, for every $i\in\mathbb N$ there exist $l_i\in\mathbb N$ such that $f^{l_i}(d_{i0})=d$. Then, for each $i\in\mathbb N$, $p\in\mathbb N^*$ $f^p(d_{i0})=d$ and $f^p(d_0)=d_1$. So, we conclude that $f^p$ is discontinuous at $d_0$ for every $p\in\mathbb N^*$.  Since there are periodic points of arbitrarily large period, the Ellis semigroup $E(\omega^3+1,f)$
is homeomorphic to $2^\mathbb N$ (by Theorem 2.7 of \cite{gru2016}).
\end{ejem}

Concerning the previous example we formulate the following question.

\begin{question} Is there a continuous function $f: \omega^3+1 \to \omega^3+1$ such  $(\omega^3+1,f)$ is transitive  and  $E(\omega^3+1,f)^*$ contains both continuous and discontinuous functions?
\end{question}

With respect to Theorem \ref{homeo}, it is natural to ask the following.

\begin{question}\label{final} Let $(\omega^\alpha+1,f)$ be a transitive dynamical systems where $\alpha \geq 1$ is a countable ordinal. If $1 < \sup\{|\omega_f(x)| : x \in (\omega^\alpha+1)' \} < \omega$,  is  $E(\omega^\alpha+1,f)$  homeomorphic to $\omega^\beta+1$ for some countable ordinal $\beta \geq 1$?
\end{question}

\bibliographystyle{amsplain}


\end{document}